\documentclass[12pt,
]{article}

\usepackage{amsmath,amsfonts,amssymb,latexsym}
\usepackage[dvips]{graphicx}
\usepackage{color}

\newtheorem{theorem}{Theorem}[section]
\newtheorem{lemma}[theorem]{Lemma}
\newtheorem{proposition}[theorem]{Proposition}

\newtheorem{definition}[theorem]{Definition}

\def\E{\mathbf{E}}
\def\P{\mathbf{P}}
\def\R{\mathbb{R}}
\def\Rp{\mathbb{R}_+}
\def\N{\mathbb{N}}

\def\sM{{\mathbf M}}
\def\sMa{{\mathbf M}_a}
\def\sMF{{\mathbf M}_{\chi}}
\def\sMFa{{\mathbf M}_0}

\def\sR{{\cal R}}

\def\sZ{{\cal Z}}
\def\sZhat{\widehat\sZ}
\def\sZstar{\sZ^*}

\def\sV{{\cal V}}

\def\Xhat{\widehat{X}}

\def\Vhat{\widehat{V}}

\def\What{\widehat{W}}

\def\sVhat{\widehat{\sV}}

\def\Zhat{\widehat{Z}}
\def\Ehat{\widehat{E}}
\def\Ebar{\bar{E}}

\def\la{\langle}
\def\ra{\rangle}
\def\wk{\stackrel{w}{\rightarrow}}
\def\dist{\stackrel{d}{=}}
\def\muz{{\bf 0}}
\def\ptime{[0,\infty)}

\def\load{{\mathcal V}}
\def\loadhat{\widehat\load}

\def\xstar{x^*}

\title{\small\bf DIFFUSION LIMITS FOR  SHORTEST REMAINING \\ PROCESSING TIME QUEUES}

\author{ {\small\sc 
H.\ Christian Gromoll\thanks{Research supported in part by NSF grant DMS
0707111}, \L ukasz Kruk, Amber L.\ Puha} \\
{\em\footnotesize University of Virginia, Maria Curie-Sk\l odowska
University,}\\ {\em\footnotesize California State University San Marcos} 
}

\mathchardef\emptyset="001F

\begin{document}

\maketitle

\begin{abstract}
We present a heavy traffic analysis for a single server queue with renewal
arrivals and generally distributed i.i.d.\ service times, in which the server
employs the Shortest Remaining Processing Time (SRPT) policy.  Under typical
heavy traffic assumptions, we prove a diffusion limit theorem for a
measure-valued state descriptor, from which we conclude a similar theorem for
the queue length process. These results allow us to make some observations on
the queue length optimality of SRPT. In particular, they provide the sharpest
illustration of the well-known tension between queue length optimality and
quality of service for this policy.  
\end{abstract}



\section{Introduction}\label{intro}
\setcounter{equation}{0}

In a single server queue employing the Shortest Remaining Processing Time
(SRPT) policy, preemptive priority is given to the job that can be completed
first, that is, the job with the shortest remaining processing time. More
precisely, consider a single server queue with renewal arrivals and i.i.d.\
service times, and let $\mathcal{I}(t)$ index in the order of their arrival
those jobs that are in the queue at time $t$.
For $i\in\mathcal{I}(t)$, let $w_i(t)$ denote the {\em residual
service time} at time $t$ of job $i$. This is the remaining amount of
processing time required to complete this job. If $j\in\mathcal{I}(t)$ is the
smallest index such that $w_j(t)\le w_i(t)$ for all $i\in\mathcal{I}(t)$, then
under SRPT, $\frac{d}{dt}w_j(t+)=-1$ and $\frac{d}{dt}w_i(t+)=0$ for all
$i\in\mathcal{I}(t)\setminus j$. 

Interest in the SRPT policy goes back to the first optimality
result of Schrage \cite{sch68}, who showed that SRPT minimizes
the number of jobs in the system, or queue length, at each point in time (see
also Smith \cite{smi76}). More explicitly, given fixed arrival and service
processes, if $Z(t)$ is the queue length at time $t$
under SRPT and $Q(t)$ is the queue length at $t$ under an arbitrary work
conserving policy, then almost surely,
\begin{equation}
  Z(t)\le Q(t),\qquad\text{for all $t\ge0$.}
  \label{minimize}
\end{equation}
This holds with no distributional assumptions on the underlying arrival and
service processes.  

Expressions for the mean response time for an M/G/1 SRPT queue were developed
earlier by Schrage and Miller \cite{schmil66}, and extended later in
Schassberger \cite{sch90} and Perera \cite{per93} (see Schreiber \cite{sch93}
for a survey of the same time period). Another notable contribution was made
by Pavlov \cite{ref:Pa} and Pechinkin \cite{ref:Pe}, who characterized the
heavy traffic limit of the steady state distributions for the queue length of
an M/G/1 SRPT queue.

Recently, there has been renewed interest in the SRPT policy, mainly in
computer science.  For example, Bansal and Harchol-Balter \cite{banhar01}
study fairness for SRPT (\cite{banhar01} is also a good source for a more
extended list of prior work on SRPT).  More recent work seeks to provide a
framework for comparing policies in the M/G/1 setting; see for example Wierman
and Harchol-Balter \cite{wiehar03}.

There has also been a recent body of work on the tail behavior of
single server queues under SRPT; see for example N\'{u}\~{n}ez Queija
\cite{nun02} and Nuyens and Zwart \cite{nuyzwa06}. They discuss the
advisability of implementing SRPT using large deviations techniques.

In \cite{dowwu06}, Down and Wu employ diffusion limits to show certain
optimality properties of a multi-layered round robin routing policy for a
system of parallel servers, each operating under SRPT. This was done under the
assumption of a finitely supported service time distribution, mainly due to
the absence at the time of diffusion limits for more general service time
distributions. In the case of a general service time distribution, Down,
Gromoll, and Puha \cite{ref:DGP} developed fluid limits for SRPT queues, and
used these to obtain a formula for state-dependent response times (on fluid
scale) of jobs entering the system (see also \cite{ref:sig}).

In this paper, we prove a diffusion limit theorem that holds for a general
service time distribution, under usual heavy traffic assumptions.  We do this
for a measure-valued state descriptor, so that diffusion limits for various
other performance measures may be obtained as corollaries; see Theorem
\ref{thrm:dlt}. In particular, we obtain a diffusion limit theorem for the
queue length process. This result reveals just how optimal SRPT is, in the sense of
\eqref{minimize}, and is explained below. 

Let $\widehat Z^r(t)=r^{-1}Z^r(r^2t)$, $t\ge0$, be the $r$th diffusion scaled
queue length process from an $r$-indexed sequence of SRPT models, as detailed
in Section \ref{sec:dl}.  In particular, we assume the fairly standard heavy
traffic assumptions \eqref{eq:exog}, \eqref{eq:service}, \eqref{eq:HTC},
\eqref{eq:finite}, \eqref{eq:ConvergenceOfInitialCondition}, and
\eqref{eq:LimitingInitialConditionFinite}.  We use $W^*(\cdot)$ to denote the
limit in distribution of the corresponding sequence of diffusion scaled
workload processes (see \eqref{eq:ConvWork}).  As noted there, $W^*(\cdot)$ is
the same for all work conserving policies and is a reflected Brownian motion
in $\mathbb{R}_+$ \cite{ref:IW}.   We use $\nu$ to denote the limiting service
time distribution (see \eqref{eq:service}) and $x^*$ to denote the supremum of
the support of $\nu$.  Informally, $\xstar$ is the largest possible job size.
Then,
\begin{theorem}\label{thrm:1}
As $r\to\infty$, the processes $\widehat Z^r(\cdot)$ converge in distribution to 
$$
Z^*(\cdot)\dist
\begin{cases}
\frac{W^*(\cdot)}{\xstar},&\hbox{if }\xstar<\infty,\\
0,&\hbox{if }\xstar=\infty.\end{cases}
$$
\end{theorem}
This result follows from Theorem \ref{thrm:dlt} by the continuous mapping
theorem. 

Theorem \ref{thrm:1} makes a striking statement about the queue length
optimality of SRPT. Consider the following simple lower bound, valid for any
work conserving policy and service time distribution $\nu$.
Assume for the moment that $x^*<\infty$. Let $Q(t)$, $t\ge0$, be the queue
length process under an arbitrary work conserving policy. Then at each time
$t\ge0$, the workload $W(t)$ is bounded above by $Q(t)x^*$, because it is the
sum of $Q(t)$ residual service times, each of which is bounded above by $x^*$.
So almost surely,
\begin{equation}
  Q(t)\ge \frac{W(t)}{x^*}, \qquad \text{for all $t\ge0$}.
  \label{static}
\end{equation}
Note that \eqref{static} makes sense when $x^*=\infty$ as well, as the right
side is interpreted as zero.  

Unlike \eqref{minimize}, which gives a universal lower bound (over all work
conserving policies) in terms of the queue length process of one such policy,
\eqref{static} gives a universal lower bound in terms of the common workload
process of all such policies. In particular, we may combine these bounds and
have, almost surely, 
\begin{equation*}
  \frac{W(t)}{x^*}\le Z(t)\le Q(t), \qquad \text{for all $t\ge0$}.
  \label{both} \end{equation*}

The bound \eqref{static} is intuitively appealing because it results from the
hypothetical configuration of residual service times that minimizes the queue
length at time $t$, given the workload at $t$. At each $t\ge0$, the queue
length minimizing configuration is the one that puts as many residual service
times as possible at $x^*$, such that they sum to $W(t)$. (To be precise,  all
of them if $x^*$ divides $W(t)$ and all but one of them otherwise).
Additionally, since the workload process is a much simpler object than the 
queue length process under SRPT, \eqref{static} may be easier to work with in
practice, when $x^*<\infty$, than \eqref{minimize}.

Of course, this bound is hypothetical because no work conserving policy,
including SRPT, can achieve such optimal configurations for all $t\ge0$,
although many may achieve it for some $t$ (including for example all times $t$
for which $W(t)=0$).  The interesting fact contained in Theorem \ref{thrm:1}
is that, on diffusion scale in heavy traffic, SRPT actually achieves the
hypothetical lower bound asymptotically, almost surely for all $t\ge0$. 

So SRPT is not only better than any other work conserving policy in the sense
of \eqref{minimize}, it is in fact as optimal as possible in the heavy traffic
limit.  Of course, this optimality is from the point of view of the server,
who one imagines wants to minimize queue length. As is well known, SRPT
performs poorly from the point of view of large jobs (see e.g.\ \cite{ref:DGP}),
who wish to minimize their time in queue, but tend to wait for long periods as
they are preempted by smaller jobs. Indeed the queue length optimality of SRPT
comes at the expense of long sojourn times for large jobs, and this tension is
made explicit by Thereom 3.1, which gives the measure-valued diffusion limit.
From this result, we see that in the heavy traffic limit, all mass is
concentrated at $x^*$. So asymptotically for all $t\ge0$, the queue consists
entirely of jobs of the largest possible size, whereas smaller jobs are
flushed out instantly.  That is, the diffusion limit in Theorem 3.1 puts the
contrast between queue length optimality and poor performance for large jobs
in the sharpest light.

In the remainder of the paper, we give a precise definition of the stochastic
model for an SRPT queue (Section \ref{sec:SRPTqueue}), state our assumptions
and main result (Section \ref{sec:dl}), and provide the proofs (Section
\ref{sec:pf}).

\subsection{Notation}\label{sec:not}

The following notation will be used throughout the paper.  Let $\N$ denote the
set of positive integers and let $\R$ denote the set of real numbers.  For
$a,b\in\R$, we write $a\vee b$ for the maximum of $a$ and $b$, and $\lfloor
a\rfloor$ for the largest integer less than or equal to $a$.  The nonnegative
real numbers $[0,\infty )$ will be denoted by $\Rp$. By convention, a sum of
the form $\sum_{i=n}^m$ with $n>m$, or a sum over an empty set of indices
equals zero.  The sets $(a,b)$, $[a,b)$, and $(a,b]$ are empty for
$a,b\in[0,\infty]$ with $a\ge b$.  For a Borel set $B\subset\Rp$, we denote
the indicator of the set $B$ by $1_{B}$.  We also define the real valued
function $\chi (x)=x$, for $x\in\Rp$.

Let $\sM$ denote the set of finite, nonnegative Borel measures on $\Rp$.  For
$\xi\in\sM$ and a Borel measurable function $g:\Rp\rightarrow\mathbb{R}$ that
is integrable with respect to $\xi$, define $\la g,\xi\ra =\int_{{\mathbb
R}_+}g(x)\xi(dx).$ The set $\sM$ is endowed with the weak topology. That is,
for  $\xi_n,\xi\in\sM$, we have $\xi_n\wk\xi$ if and only if $\la g,\xi_n\ra
\rightarrow \la g,\xi\ra$ as $n\rightarrow\infty$, for all
$g:\Rp\rightarrow\R$ that are bounded and continuous. With this topology,
$\sM$ is a Polish space \cite{ref:P}.    We denote the zero measure in $\sM$
by $\muz$ and the measure in $\sM$ that puts one unit of mass at the point
$x\in\Rp$ by $\delta_x$.  For $x\in\Rp$, the measure $\delta_x^+$ is
$\delta_x$ if $x>0$ and $\muz$ otherwise.  For $\xi\in\sM$, we say that
$x\in\Rp$ is a $\xi$-continuity point if $\langle 1_{\{x\}},\xi\rangle=0$.
Let $\sMa$ denote those elements of $\sM$ that do not charge the origin.  We
say that a measure $\xi\in\sM$ has a finite first moment if
$\la\chi,\xi\ra<\infty$.  Let $\sMF$ denote the set of all such measures and
let $\sMFa=\sMF\cap\sMa$.

We use ``$\overset{d}{=}$'' for equality in distribution and ``$\Rightarrow$''
to denote convergence in distribution of random elements of a metric space.
Unless otherwise specified, all stochastic processes used in this paper are
assumed to have paths that are right continuous with finite left limits
(r.c.l.l.). For a Polish space ${\cal S}$, we denote by
$\mathbf{D}([0,\infty),{\cal S})$ the space of r.c.l.l.\ functions from
$[0,\infty)$ into ${\cal S}$, endowed with the Skorohod $J_1$-topology
\cite{ref:EK}.

\section{Stochastic Model for an SRPT Queue}\label{sec:stochModel}\label{sec:SRPTqueue}
\setcounter{equation}{0}
Our stochastic model of an SRPT queue consists of the following:
a random initial
condition $\sZ(0)\in\sM$ specifying the state of the system at time zero,
stochastic primitives $E(\cdot)$ and $\{v_k\}_{k\in\N}$ describing the
arrival of jobs to the queue and their service times,
and a measure valued state descriptor $\sZ(\cdot)$ describing the time
evolution of the system.  These are defined below.

\paragraph{Initial condition.} The initial condition specifies the number
$Z(0)$ of jobs in the queue at time zero, as well as the initial service
time of each job.  Assume that $Z(0)$ is a nonnegative integer valued
random variable that is finite almost surely.  The initial service times
are the first $Z(0)$ elements of a sequence $\{\tilde v_j\}_{j\in\N}$
of strictly positive, finite random variables.  The initial job with
service time $\tilde v_j$, $j\le Z(0)$, is called job $j$.

A convenient way to express the initial condition is to define an initial
random measure $\sZ(0)\in\sM$ by
$$
\sZ(0)=\sum_{j=1}^{Z(0)} \delta_{\tilde v_j},
$$
which equals $\mathbf{0}$ if $Z(0)=0$.  Our assumptions imply that
$\sZ(0)$ satisfies
\begin{equation}\label{eq:InitialConditionFinite}
\mathbf{P}(\la 1,\sZ(0)\ra\vee\la \chi,\sZ(0)\ra<\infty)=1.
\end{equation}
In particular, the number of initial jobs and the initial workload are
finite almost surely, and so $\sZ(0)\in\sMFa$ almost surely.

\paragraph{Stochastic primitives.}
The stochastic primitives consist of an exogenous arrival process
$E(\cdot)$ and a sequence of initial service times $\{v_k\}_{k\in\N}$.
The arrival process $E(\cdot)$ is a rate $\alpha\in(0,\infty)$ delayed
renewal process such that the interarrival times have standard deviation $a\in[0,\infty)$.
For $t\in\ptime$, $E(t)$ represents the number of jobs
that arrive to the queue during the time interval $(0,t]$.  Jobs arriving
after time $0$ are indexed by integers $j>Z(0)$.  For $t\in\ptime$, let
\begin{equation}\label{def:A}
A(t)=Z(0)+E(t).
\end{equation}
Then job $j\in\N$ arrives at time $T_j=\inf\{t\in\ptime: A(t)\ge j\}$.  Hence,
for  $i<j$, $T_i\le T_j$ and we say that job $i$ arrives before job $j$.

For each $k\in\N$, the random variable $v_k$ represents the initial
service time of the $(Z(0)+k)$th job. That is, job $j>Z(0)$ has
initial service time $v_{j-Z(0)}$. Assume that the random variables
$\{v_k\}_{k\in\N}$ are strictly positive and form an independent
and identically distributed sequence with common Borel distribution
$\nu$ on $\Rp$.  Assume that the mean $\la\chi,\nu\ra\in(0,\infty)$
and standard deviation $b=\sqrt{\la\chi^2,\nu\ra -\la\chi,\nu\ra^2}\in[0,\infty)$.
Let $\beta=\la\chi,\nu\ra^{-1}$.  Define the traffic intensity 
$\rho=\alpha/\beta$.

It will be convenient to combine the stochastic primitives into a single,
measure valued load process.
\begin{definition}\label{d.loadprocess}  The load process is given by
$$
\load(t)=\sum_{k=1}^{E(t)}\delta_{v_k}, \quad \text{for $t\in\ptime$}.
$$
\end{definition}
Then $\load(\cdot)\in \mathbf{D}(\ptime,\sM)$ since $E(\cdot)\in\mathbf{D}(\ptime,\Rp)$.

\paragraph{Evolution of the residual service times.}
In an SRPT queue, the smallest nonzero residual service time decreases at
rate one until either it becomes zero or a job arrives that has a smaller
initial service time, at which time the rate changes to zero and the new
smallest nonzero residual service time begins decreasing at rate one.
We adopt the convention that in case of a tie, the residual service
time of the job that arrived first (that is, the job with smaller index)
begins decreasing at rate one. 

For $j\in\N$ and $t\in\ptime$, let $w_j(t)$ denote the residual service time
of job $j$.  By convention, for $j\in\N$ and $t\in[0,T_j]$,
$$
w_j(t)=\begin{cases} \tilde v_j,&1\le j\le Z(0),\\ v_{j-Z(0)},&j>Z(0).\end{cases}
$$
Furthermore, for $j\in\N$, if $D_j$ denotes the time at which job $j$
completes service and departs the system, then $w_j(t)=0$ for all $t\ge D_j$.
On $(T_j,D_j)$, $w_j(\cdot)$ is nonincreasing.  In particular, $w_j(\cdot)$
decreases at rate one when job $j$ is in service, and is constant when job $j$
is not in service.   See \cite{ref:DGP} for a detailed definition of the
residual service times.

\paragraph{Measure-valued state descriptor.}
For $t\in\ptime$, define the state descriptor by
\begin{equation}\label{eq:z}
\sZ(t)=\sum_{j=1}^{A(t)}\delta^+_{w_j(t)}.
\end{equation}

\section{Diffusion Limit Theorem}\label{sec:dl}
\setcounter{equation}{0}

We first define a sequence of systems over which the limit
is taken.  Let $\sR$ be a sequence of positive real numbers
increasing to infinity.  Consider an $\sR$-indexed sequence
of stochastic models, each defined as in Section \ref{sec:stochModel}.
For each $r\in\sR$, there is an initial condition ${\mathcal
Z}^r(0)$; there are stochastic primitives $E^r(\cdot)$ and
$\{v_k^r\}_{k\in\N}$ with parameters $\alpha^r$, $a^r$, $\nu^r$, $\beta^r$, $b^r$, and
$\rho^r$, and an arrival process $A^r(\cdot)$ with arrival times
$\{T_j^r\}_{j\in\N}$; there is a corresponding measure valued load process
$\load^r(\cdot)$;  there is a state descriptor ${\mathcal Z}^r(\cdot)$.
The stochastic elements of each model are defined on a probability space
$(\Omega^r,{\mathcal F}^r,\P^r)$ with expectation operator $\E^r$.
A diffusion scaling (or central limit theorem scaling) is applied to each
model in the $\sR$-indexed sequence as follows.  For each
$r\in\sR$ and $t\in\ptime$, let
\begin{equation}\label{e.firstflscaling}
\Ehat^r(t)=\frac{1}{r}\left(E^r(r^2t)-r^2t\alpha^r\right).
\end{equation}
Also, for each $r\in\sR$ and $t\in\ptime$, let
\begin{equation}\label{e.firstflscaling'}
\sZhat^r(t)=\frac{1}{r}\sZ^r(r^2t)
\qquad\hbox{and}\qquad
\What^r(t)=\la \chi,\sZhat^r(t)\ra.
\end{equation}

Let $\alpha, a\in(0,\infty)$ and define $\alpha(t)=\alpha t$ for all
$t\in\ptime$.  Let $\nu$ be a probability measure such that
\begin{equation}\label{eq:servicelimit}
\la 1_{\{0\}},\nu\ra=0,\qquad
\la \chi,\nu\ra=1/\alpha,
\qquad\hbox{and}\qquad
0<\la \chi^2,\nu\ra<\infty.
\end{equation}
Set $b=\sqrt{ \la \chi^2,\nu\ra-\la\chi,\nu\ra^2}$ and
$$
\xstar=\sup\{ x\in\Rp : \la 1_{[0,x]},\nu\ra <1\}.
$$

For the sequence of stochastic primitives we make the following
asymptotic assumptions.  For the exogenous arrival processes,
assume that as $r\to\infty$,
\begin{equation}\label{eq:exog}
\alpha^r\to\alpha,\qquad a^r\to a,
\qquad\hbox{and}\qquad
\Ehat^r(\cdot)\Rightarrow E^*(\cdot),
\end{equation}
where $E^*(\cdot)$ is a Brownian motion starting from zero with drift zero and
variance $a^2\alpha^3$ per unit time.
This implies a functional weak law of large numbers for the exogenous arrival
processes.  In particular, it implies that as $r\to\infty$,
$$
\Ebar^r(\cdot)\Rightarrow \alpha(\cdot),
$$
where $\Ebar^r(t)=E^r(r^2t)/r^2$ for all $t\in\ptime$ and $r\in\sR$.
For the sequence of service time distributions, assume that as $r\to\infty$,
\begin{equation}\label{eq:service}
\nu^r\wk\nu
\qquad\hbox{and}\qquad
\la \chi^2,\nu^r\ra\to\la \chi^2,\nu\ra.
\end{equation}
Then $\beta^r\to\alpha$, $\rho^r\to 1$,
and $b^r\to b$ as $r\to\infty$. It also follows that 
$\{\nu^r, r\in\sR\}$ satisfies a Lindeberg-Feller condition, i.e., for all $\varepsilon>0$,
\begin{equation}\label{eq:LF}
\lim_{r\to\infty}\la (\chi-\la\chi,\nu^r\ra)^2\left(1_{[0,\la\chi,\nu^r\ra-\varepsilon r)}+1_{(\la\chi,\nu^r\ra+\varepsilon r,\infty)}\right),\nu^r\ra=0.
\end{equation}
In addition, assume the heavy traffic condition that for some $\gamma\in\R$,
\begin{equation}\label{eq:HTC}
\lim_{r\to\infty} r(1-\rho^r)=\gamma.
\end{equation}
Finally, if $\xstar<\infty$, also assume that for all $x>\xstar$,
\begin{equation}\label{eq:finite}
\lim_{r\to\infty} r\la \chi1_{(x,\infty)},\nu^r\ra=0.
\end{equation}

For the sequence of diffusion scaled initial conditions $\{\sZhat^r(0) : r>0\}$,
assume that as $r\to\infty$,
\begin{equation}\label{eq:ConvergenceOfInitialCondition}
\What^r(0)\Rightarrow W_0^*,
\end{equation}
for some random variable $W_0^*$.  Then from \eqref{eq:exog}, \eqref{eq:service} (which implies \eqref{eq:LF}), \eqref{eq:HTC}, \eqref{eq:ConvergenceOfInitialCondition}, and the fact that SRPT is a work conserving discipline, it follows that, as $r\to\infty$,
\begin{equation}\label{eq:ConvWork}
\What^r(\cdot)\Rightarrow W^*(\cdot),
\end{equation}
where $W^*(\cdot)$ is a reflected Brownian motion with initial value
$W^*(0)\dist W_0^*$, variance $(a^2+b^2)\alpha$ per unit time, and drift $-\gamma$
(see \cite{ref:IW}).
Further assume that, as $r\to\infty$,
\begin{equation}\label{eq:LimitingInitialConditionFinite}
\sZhat^r(0) \Rightarrow
\begin{cases}
\frac{W_0^*}{\xstar}\delta_{\xstar},&\hbox{if }\xstar<\infty,\\
\muz,&\hbox{if }\xstar=\infty.
\end{cases}
\end{equation}
Note that \eqref{eq:LimitingInitialConditionFinite} implies that $\sZhat^r(0)$ converges in distribution to a random measure
that is almost surely an invariant state (see \cite[Corollary 3.7]{ref:DGP}).

\begin{theorem}\label{thrm:dlt}  Under the asymptotic assumptions
\eqref{eq:exog}, \eqref{eq:service}, \eqref{eq:HTC}, \eqref{eq:finite}, \eqref{eq:ConvergenceOfInitialCondition}, and \eqref{eq:LimitingInitialConditionFinite}, the
sequence $\{\sZhat^r(\cdot):r\in\sR\}$ converges in distribution
on $\mathbf{D}([0,\infty),\sM)$ to a measure valued process $\sZstar(\cdot)$ such
that
$$
\sZstar(\cdot)\dist
\begin{cases}
\frac{W^*(\cdot)}{\xstar}\delta_{\xstar},&\hbox{if }\xstar<\infty,\\
\muz,&\hbox{if }\xstar=\infty.\end{cases}
$$
\end{theorem}

This result, in the first case when $\xstar<\infty$, is a continuous analog of
the diffusion limit result for a multi-class static buffer priority queue,
where in the diffusion limit work only resides in the lowest priority class
\cite{ref:W2}.  In an SRPT queue, those jobs with larger service times receive
lower priority.  Hence, an informal restatement of the first case is that in
the diffusion limit the work concentrates in jobs with the largest possible
service time, i.e., the lowest priority.  The case when $\xstar=\infty$ is the
natural extension of this result when there is no largest possible service
time.  Indeed, for the work to get pushed out to infinity on diffusion scale
while the diffusion scaled workload process converges, the queue length must
necessarily tend to zero.

\section{Proofs}\label{sec:pf}
\setcounter{equation}{0}

Throughout this section we assume that \eqref{eq:exog}, \eqref{eq:service},
\eqref{eq:HTC}, \eqref{eq:finite}, \eqref{eq:ConvergenceOfInitialCondition},
and \eqref{eq:LimitingInitialConditionFinite} hold.  In Section
\ref{sec:basic}, we state a well known result  and use it to derive three
diffusion limit results to be used in the sequel.  In Section \ref{sec:main},
Theorem \ref{thrm:dlt} is proved.

\subsection{Diffusion Limits for Load Related Processes}\label{sec:basic}

The following result is well known and follows from \cite[Theorem 3.1]{ref:P} used to extend \cite[Section 17.3]{ref:B}.

\begin{proposition}\label{lem:basic}
For each $r\in\sR$, let $\{x_k^r\}_{k=1}^{\infty}$ be an independent and
identically distributed sequence of nonnegative random variables on
$(\Omega^r, \mathcal{F}^r,\mathbf{P}^r)$ with finite mean $\mu^r$ and standard
deviation $\sigma^r$, that is independent of $E^r(\cdot)$.  Suppose that for
some finite nonnegative constants $\mu$ and $\sigma$, $\mu^r\to\mu$ and
$\sigma^r\to\sigma$ as $r\to\infty$.  Further assume that for each
$\varepsilon>0$,
$$
\lim_{r\to\infty}\E^r\left[ (x_1^r-\mu^r)^2 ; \left| x_1^r-\mu^r\right| > r\varepsilon \right]=0.
$$
For $r\in\sR$, $n\in\N$, and $t\in\ptime$, let
$$
X^r(n)=\sum_{k=1}^n x_k^r.
\qquad\hbox{and}\qquad
\Xhat^r(t)=\frac{X^r(\lfloor r^2t\rfloor ) -\lfloor r^2t\rfloor \mu^r}{r}.
$$
Then as $r\to\infty$, $(\Ehat^r(\cdot),\Xhat^r(\cdot))\Rightarrow (E^*(\cdot),X^*(\cdot))$,
where $E^*(\cdot)$ is given by \eqref{eq:exog} and $X^*(\cdot)$ is a Brownian motion starting from zero with zero drift and
variance $\sigma^2$ per unit time, that is independent of $E^*(\cdot)$.
Furthermore, as $r\to\infty$,
$$
\frac{X^r(r^2\Ebar^r(\cdot))-r^2\alpha^r(\cdot)\mu^r}{r}
\Rightarrow
X^*(\alpha(\cdot))+\mu E^*(\cdot),
$$
where for each $r\in\sR$ and $t\in\ptime$, $\alpha^r(t)=\alpha^r t$.
\end{proposition}

Note that the limiting process $X^*(\alpha(\cdot))+\mu E^*(\cdot)$ in
Proposition \ref{lem:basic} is a Brownian motion starting from zero with zero
drift and variance $\alpha\sigma^2+\mu^2\alpha^3 a^2$ per unit time.  We apply this
proposition to three processes of interest here, that we respectively refer to
as the total load, the truncated load, and the tail load processes.  For 
$r\in\sR$ and $t\in\ptime$, let
$$
\loadhat^r(t)=\frac{1}{r}\left(\load^r(r^2t)-r^2\alpha^r t\nu^r\right).
$$
Then, for $r\in\sR$, let the total load and scaled total load processes
be given respectively by 
$$
V^r(\cdot)=\la\chi,\load^r(\cdot)\ra\qquad\hbox{and}\qquad\Vhat^r(\cdot)=\la\chi,\loadhat^r(\cdot)\ra.
$$
Then, for $r\in\sR$,
$$
\Vhat ^r(\cdot)=\frac{\sum_{k=1}^{r^2\Ebar^r(\cdot)} v_k^r -r^2\alpha^r(\cdot)\la \chi,\nu^r\ra}{r}.
$$
From \eqref{eq:service} and Proposition \ref{lem:basic}, it follows that as $r\to\infty$,
$$
\Vhat^r(\cdot)\Rightarrow V^*(\cdot),
$$
where $V^*(\cdot)$ is a Brownian motion starting from zero with zero drift
and variance $\alpha(a^ 2 +b^2)$ per unit time.

Next we consider the truncated load process.
For $r\in\sR$ and $x\in\Rp$,
let
$$
V_x^r(\cdot)=\la \chi 1_{[0,x]},\sV^r(\cdot) \ra
\qquad\hbox{and}\qquad
\Vhat_x^r(\cdot)=\la \chi 1_{[0,x]}, \sVhat^r(\cdot)\ra.
$$
Then, for $r\in\sR$ and $x\in\Rp$,
$$
\Vhat_x ^r(\cdot)=\frac{\sum_{k=1}^{r^2\Ebar^r(\cdot)} v_k^r1_{\{v_k^r\le x\}} -r^2\alpha^r(\cdot)\la \chi1_{[0,x]},\nu^r\ra}{r}.
$$
Note that \eqref{eq:service} implies that for any 
$\nu$-continuity point $x\in\Rp$, as $r\to\infty$,
\begin{equation}\label{eq:TruncService}
\la \chi^2 1_{[0,x]},\nu^r\ra\to\la \chi^2 1_{[0,x]},\nu\ra.
\end{equation}
Hence \eqref{eq:service} and Proposition \ref{lem:basic} imply that for any
$\nu$-continuity point $x\in\Rp$, as $r\to\infty$,
\begin{equation}\label{eq:TruncDiff}
\Vhat_x^r(\cdot)\Rightarrow V_x^*(\cdot),
\end{equation}
where $V_x^*(\cdot)$ is a Brownian motion starting from zero with drift zero
and finite variance per unit time.

Finally, for each $r\in\sR$ and $x\in\Rp$, we consider the tail load process
$V^r(\cdot)-V_x^r(\cdot)$.  Then, for $r\in\sR$ and $x\in\Rp$,
$$
\Vhat^r(\cdot)-\Vhat_x ^r(\cdot)=\frac{\sum_{k=1}^{r^2\Ebar^r(\cdot)} v_k^r1_{\{v_k^r>x\}} -r^2\alpha^r(\cdot)\la \chi1_{(x,\infty)},\nu^r\ra}{r}.
$$
Note that \eqref{eq:service} (which implies \eqref{eq:TruncService}) also
implies that for any $\nu$-continuity point $x\in\Rp$, as $r\to\infty$,
$$
\la \chi^2 1_{(x,\infty)},\nu^r\ra\to\la \chi^2 1_{(x,\infty)},\nu\ra.
$$
Hence, \eqref{eq:service} and Proposition \ref{lem:basic} imply that for any
$\nu$-continuity point $x\in\Rp$, as $r\to\infty$,
\begin{equation}\label{eq:TailDiff}
\Vhat^r(\cdot)-\Vhat_x^r(\cdot)\Rightarrow T^*_x(\cdot),
\end{equation}
where $T^*_x(\cdot)$ is a Brownian motion starting from zero with drift zero and variance $s_x^2$ per unit time.
Here,
\begin{equation}\label{eq:var}
s_x^2
=\alpha(\la \chi^2 1_{(x,\infty)},\nu\ra - \la \chi 1_{(x,\infty)},\nu\ra^2)
+ \la \chi 1_{(x,\infty)},\nu\ra^2 
\alpha^3 a^2.
\end{equation}
Notice that if $\xstar<\infty$ and $x>\xstar$,
then $x$ is a $\nu$-continuity point and $\la 1_{(x,\infty)},\nu\ra=0$.
Hence, if $x>\xstar$, then in \eqref{eq:var}, $s_x^2=0$, i.e.,
\begin{equation}\label{eq:T}
T^*_x(\cdot)\equiv 0.
\end{equation}

\subsection{Proof of the Main Theorem}\label{sec:main}

Here we use the diffusion limits for the load related processes
derived in Section \ref{sec:basic} to prove the main result.  We use the
result about the scaled truncated load process to prove that, on diffusion
scale, the truncated queue length tends to zero when the truncation is below
$x^*$, the supremum of the support of the limiting service time distribution.
Then we use the result about the scaled tail load processes to prove that, on
diffusion scale, the queue length above $x$ tends to zero when $x$ is above
$x^*$. Then these two results are put together to show that in the diffusion
limit, the queue mass concentrates at $x^*$.

For $r\in\sR$ and $x\in\Rp$, let
\begin{eqnarray}
Z_x^r(\cdot)=\la 1_{[0,x]},\sZ^r(\cdot)\ra
&\qquad\hbox{and}\qquad&
W_x^r(\cdot)=\la \chi 1_{[0,x]},\sZ^r(\cdot)\ra,\label{def:Wx}\\
\Zhat_x^r(\cdot)=\la 1_{[0,x]},\sZhat^r(\cdot)\ra
&\qquad\hbox{and}\qquad&
\What_x^r(\cdot)=\la \chi 1_{[0,x]},\sZhat^r(\cdot)\ra.\label{def:hatWx}
\end{eqnarray}

\begin{lemma}\label{lem:ZeroQueue}
For any $x\in(0,\xstar)$, as $r\to\infty$,
\begin{equation}\label{eq:MassToZero}
\Zhat_x^r(\cdot) \Rightarrow 0.
\end{equation}
\end{lemma}

\noindent{\bf Proof}.  Since $\Zhat_y^r(\cdot)\le \Zhat_x^r(\cdot)$
for each $0< y\le x<\xstar$, it suffices to verify
\eqref{eq:MassToZero} for $x\in(0,\xstar)$ that are $\nu$-continuity
points.  Fix such an $x$.  For $r\in\sR$ and $t\in\ptime$,
let
$$
\tau_x^r(t)=\sup\{ s\in[0,t] : \Zhat_x^r(s)=0\},
$$
which is taken to be zero if $\{ s\in[0,t] : \Zhat_x^r(s)=0\}=\emptyset$.
Then, for $r\in\sR$ and $t\in\ptime$,
\begin{equation}\label{eq:KeyINEQ1}
\begin{aligned}
\Zhat_x^r(t)
&\le\Zhat_x^r(\tau_x^r(t)) + \frac{E^r(r^2t)-E^r(r^2\tau_x^r(t))}{r}\\
&=\Zhat_x^r(\tau_x^r(t)) + \Ehat^r(t)-\Ehat^r(\tau_x^r(t))
+r\left(t-\tau_x^r(t)\right)\alpha^r.
\end{aligned}
\end{equation}

First, we obtain an upper bound on
$\Zhat_x^r(\tau_x^r(\cdot))$.  Fix $r\in\sR$ and $t\in\ptime$.  Either $\tau_x^r(t)=0$ or
$\tau_x^r(t)>0$.  If $\tau_x^r(t)=0$, then $\Zhat_x^r(\tau_x^r(t))=\Zhat_x^r(0)$.
Otherwise, $\tau_x^r(t)>0$.  If $\Zhat_x^r(\tau_x^r(t))=0$, then any nonnegative upper
bound suffices.  Hence, without loss of generality, we also assume that $\Zhat_x^r(\tau_x^r(t))>0$.
Then $\Zhat_x^r(\tau_x^r(t)-)=0$ and $\Zhat_x^r(\tau_x^r(t))>0$.  Hence,
in the $r$th system at time $r^2\tau_x^r(t)$, the exogenous arrival process jumps and at least one of the
entering jobs has an initial service time in $[0,x]$, and/or the residual service
time of the job in service just before time $r^2\tau_x^r(t)$ decreases to $x$.
Therefore, $\Zhat_x^r(\tau_x^r(t))\le \Ehat^r(\tau_x^r(t))-\Ehat^r(\tau_x^r(t)-)+\frac{1}{r}$.
Combining the bounds for $\tau_x^r(t)=0$ or $\tau_x^r(t)>0$ gives
$$
\Zhat_x^r(\tau_x^r(t))\le
\Zhat_x^r(0)+\Ehat^r(\tau_x^r(t))-\Ehat^r(\tau_x^r(t)-)+\frac{1}{r},
$$
where we adopt the convention $\Ehat^r(0-) =\Ehat^r(0) =0$.
Hence, for $r\in\sR$ and $t\in\ptime$,
\begin{equation}\label{eq:KeyINEQ3}
\Zhat_x^r(t)
\le
\Zhat_x^r(0) + \Ehat^r(t)-\Ehat^r(\tau_x^r(t)-) + \frac{1}{r}
+r\left(t-\tau_x^r(t)\right)\alpha^r.
\end{equation}

For $r\in\sR$ and $t\in\ptime$, let $\theta_x^r(t)=t-\tau_x^r(t)$.
In order to show that the upper bound in \eqref{eq:KeyINEQ3}
tends to zero and thereby prove \eqref{eq:MassToZero}, 
it suffices to show that as $r\to\infty$,
\begin{equation}\label{eq:BusyZero}
r\theta_x^r(\cdot)\Rightarrow 0.
\end{equation}
To see this, assume that \eqref{eq:BusyZero} holds.  Then, for $r\in\sR$ and $t\in\ptime$, let
$$
\tilde\theta_x^r(t)=\theta_x^r(t)+\frac{1}{r^2}.
$$
By \eqref{eq:BusyZero}, as $r\to\infty$,
\begin{equation}\label{eq:TimeDiffToZero}
\theta_x^r(\cdot)\Rightarrow 0
\qquad\hbox{and}\qquad
\tilde\theta_x^r(\cdot)\Rightarrow 0.
\end{equation}
We have that for each $r\in\sR$ and $t\in\ptime$,
$$
\Ehat^r(t)-\Ehat^r(\tau_x^r(t)-)
=
\Ehat^r(t)-\frac{1}{r}E^r(r^2\tau_x^r(t)-)+r\tau_x^r(t)\alpha^r.
$$
Hence, for each $r\in\sR$ and $t\in\ptime$,
$$
\Ehat^r(t)-\Ehat^r\left(\tau_x^r(t)\right)
\le 
\Ehat^r(t)-\Ehat^r(\tau_x^r(t)-)
\le 
\Ehat^r(t)-\Ehat^r\left(\tau_x^r(t)-\frac{1}{r^2}\right)+\frac{\alpha^r}{r},
$$
where we adopt the convention that $E^r( t )=E^r(0)$ if $t<0$.
Therefore, for each $r\in\sR$ and $t\in\ptime$,
$$
\Ehat^r(t)-\Ehat^r\left(t-\theta_x^r(t)\right)
\le 
\Ehat^r(t)-\Ehat^r(\tau_x^r(t)-)
\le 
\Ehat^r(t)-\Ehat^r\left(t-\tilde\theta_x^r(t)\right)+\frac{\alpha^r}{r}.
$$
By \eqref{eq:exog}, the fact that $E^*(\cdot)$ is continuous almost surely, and
\eqref{eq:TimeDiffToZero},
it follows that, as $r\to\infty$,
$$
\Ehat^r(\cdot)-\Ehat^r\left(\cdot-\theta_x^r(\cdot)\right)\Rightarrow 0
\qquad\hbox{and}\qquad
\Ehat^r(\cdot)-\Ehat^r\left(\cdot-\tilde\theta_x^r(\cdot)\right)+\frac{\alpha^r}{r}\Rightarrow 0.
$$
(see \cite[Section 17]{ref:B}).
Hence, as $r\to\infty$,
\begin{equation}\label{eq:DiffThm}
\Ehat^r(\cdot)-\Ehat^r(\tau_x^r(\cdot)-)\Rightarrow 0.
\end{equation}
Then \eqref{eq:KeyINEQ3}, \eqref{eq:LimitingInitialConditionFinite}, \eqref{eq:DiffThm}, 
\eqref{eq:exog},  and \eqref{eq:BusyZero} together imply \eqref{eq:MassToZero}.

Hence, all that remains is to prove \eqref{eq:BusyZero}.
For this, for each $r\in\sR$ and $t\in\ptime$,
we exploit the behavior of $W_x^r(\cdot)$ (defined in \eqref{def:Wx})
on time intervals of the form $(r^2\tau_x^r(t),r^2t]$ to derive
an expression that relates $\What_x^r(t)$ and $\theta_x^r(t)$.
In particular,
since for each $r\in\sR$ and $t\in\ptime$, $Z_x^r(s)\not=0$ for all $s\in(r^2\tau^r_x(t),r^2t]$
and the service discipline is SRPT, it follows that for each $r\in\sR$ and $t\in\ptime$,
$$
W_x^r(r^2t)=W_x^r(r^2\tau_x^r(t)) +V_x^r(r^2t)-V_x^r(r^2\tau_x^r(t))-r^2(t-\tau_x^r(t)).
$$
Then, for $r\in\sR$ and $t\in\ptime$,
$$
\What_x^r(t)=\What_x^r(\tau_x^r(t)) +\Vhat_x^r(t)-\Vhat_x^r(\tau_x^r(t))
+
\left(\alpha^r\la \chi 1_{[0,x]},\nu^r\ra-1\right)
r\theta_x^r(t).
$$
Using the same line of reasoning that gave rise to \eqref{eq:KeyINEQ3},
for $r\in\sR$ and $t\in\ptime$,
$$
\What_x^r(t)\le \What_x^r(0) +\Vhat_x^r(t)-\Vhat_x^r(\tau_x^r(t)-) +\frac{x}{r}
+
\left(\alpha^r\la \chi 1_{[0,x]},\nu^r\ra-1\right)
r\theta_x^r(t).
$$
Since $\What_x^r(t)\ge 0$ for all $r\in\sR$ and $t\in\ptime$, it follows that
for $r\in\sR$ and $t\in\ptime$,
\begin{equation}\label{eq:KeyINEQ2}
\left(1-\alpha^r\la \chi 1_{[0,x]},\nu^r\ra\right)\theta_x^r(t)
\le
\frac{\What_x^r(0)}{r} +\frac{\Vhat_x^r(t)}{r}-\frac{\Vhat_x^r(\tau_x^r(t)-)}{r}+\frac{x}{r^2}.
\end{equation}
By \eqref{eq:service} and the fact that $x$ is a $\nu$-continuity point, we have that
\begin{equation}\label{eq:TrunMean}
\lim_{r\to\infty}\left(1-\alpha^r\la \chi 1_{[0,x]},\nu^r\ra\right)=1-\alpha\la \chi 1_{[0,x]},\nu\ra>0.
\end{equation}
Hence, for $r$ sufficiently large,
$\left(1-\alpha^r\la \chi 1_{[0,x]},\nu^r\ra\right)\theta_x^r(t)\ge 0$ for all $t\in\ptime$.
Then \eqref{eq:KeyINEQ2}, \eqref{eq:LimitingInitialConditionFinite}, \eqref{eq:TruncDiff},
and \eqref{eq:TrunMean} together imply that as $r\to\infty$,
$$
\theta_x^r(\cdot)\Rightarrow 0.
$$
Hence, by \eqref{eq:TruncDiff} and the same line of reasoning that gave rise
to \eqref{eq:DiffThm}, as $r\to\infty$,
$$
\Vhat_x^r(\cdot)-\Vhat_x^r(\tau_x^r(\cdot)-)\Rightarrow 0.
$$
Therefore, if one multiplies \eqref{eq:KeyINEQ2} by $r$ and
uses this and \eqref{eq:LimitingInitialConditionFinite}, \eqref{eq:BusyZero} follows.
$\hfil\Box$

We are ready to use Lemma \ref{lem:ZeroQueue}, \eqref{eq:TailDiff}, and \eqref{eq:T}
to prove the main theorem.

\noindent{\bf Proof of Theorem \ref{thrm:dlt}}. 
First suppose that $\xstar=\infty$.  Then
it suffices to show that as $r\to\infty$,
\begin{equation}\label{eq:ZeroDiff}
\Zhat^r(\cdot)\Rightarrow 0.
\end{equation}
For $r\in\sR$, $x\in\Rp$ and $t\in\ptime$, we have
\begin{align}
\Zhat^r(t)&=\Zhat_x^r(t)+\la 1_{(x,\infty)} ,\sZhat^r(t)\ra\nonumber\\
&\le
\Zhat_x^r(t)+\frac{1}{x}\la \chi 1_{(x,\infty)} ,\sZhat^r(t)\ra\nonumber\\
&\le
\Zhat_x^r(t)+\frac{1}{x}\What^r(t).\nonumber
\end{align}
Hence \eqref{eq:ZeroDiff} follows from Lemma \ref{lem:ZeroQueue}, \eqref{eq:ConvWork}, and the fact that $x$ is arbitrary.

Next suppose that $\xstar<\infty$ and let $\varepsilon>0$ be such that
$\xstar-\varepsilon$ is a $\nu$-continuity point.
Then by Lemma \ref{lem:ZeroQueue}, as $r\to\infty$,
\begin{equation}\label{eq:NoQueueBelow}
\Zhat_{\xstar-\varepsilon}^r(\cdot)
\Rightarrow
0.
\end{equation}
For $r\in\sR$ and $t\in\ptime$, we have
$$
\la 1_{(\xstar+\varepsilon,\infty)},\sZhat^r(t)\ra
\le
\frac{1}{\xstar+\varepsilon}\la \chi 1_{(\xstar+\varepsilon,\infty)},\sZhat^r(t)\ra.
$$
But, for $r\in\sR$ and $t\in\ptime$,
\begin{align*}
\la \chi 1_{(\xstar+\varepsilon,\infty)},\sZhat^r(t)\ra
&\le
\la \chi 1_{(\xstar+\varepsilon,\infty)},\sZhat^r(0)\ra
+
\frac{V^r(r^2t)-V_{\xstar+\varepsilon}^r(r^2t)}{r}\\
&\le
\la \chi 1_{(\xstar+\varepsilon,\infty)},\sZhat^r(0)\ra
 +
\Vhat^r(t)-\Vhat_{\xstar+\varepsilon}^r(t)\\
&\qquad +rt\alpha^r \la \chi 1_{(\xstar+\varepsilon,\infty)},\nu^r\ra.
\end{align*}
 Hence, by \eqref{eq:ConvergenceOfInitialCondition}, \eqref{eq:LimitingInitialConditionFinite}, \eqref{eq:TailDiff}, \eqref{eq:T}, and
 \eqref{eq:finite}, as $r\to\infty$,
\begin{equation}\label{eq:NoWorkAbove}
\la \chi 1_{(\xstar+\varepsilon,\infty)},\sZhat^r(\cdot)\ra
\Rightarrow 0.
\end{equation}
Therefore,
\begin{equation}\label{eq:NoQueueAbove}
\la 1_{(\xstar+\varepsilon,\infty)},\sZhat^r(\cdot)\ra
\Rightarrow 0.
\end{equation}
In addition, \eqref{eq:NoWorkAbove} together with \eqref{eq:NoQueueBelow}
and \eqref{eq:ConvWork} implies that as $r\to\infty$, 
\begin{equation}\label{eq:Concentrated}
\la \chi 1_{(\xstar-\varepsilon,\xstar+\varepsilon]},\sZhat^r(\cdot)\ra
\Rightarrow W^*(\cdot).
\end{equation}
Since for $r\in\sR$,
$$
\frac{1}{\xstar+\varepsilon}
\la \chi 1_{(\xstar-\varepsilon,\xstar+\varepsilon]},\sZhat^r(\cdot)\ra
\le
\la 1_{(\xstar-\varepsilon,\xstar+\varepsilon]},\sZhat^r(\cdot)\ra
\le
\frac{1}{\xstar-\varepsilon}
\la \chi 1_{(\xstar-\varepsilon,\xstar+\varepsilon]},\sZhat^r(\cdot)\ra,
$$
\eqref{eq:Concentrated}, \eqref{eq:NoQueueBelow}, \eqref{eq:NoQueueAbove}, and the fact that $\varepsilon>0$ can be made arbitrarily small completes the proof.
\hfill$\Box$

\bibliographystyle{hcg}
\bibliography{GKPFinal}

\vspace{2ex}
\begin{minipage}{2in}
\footnotesize {\sc Department of Mathematics\\
University of Virginia\\
Charlottesville, VA 22903\\
E-mail:} gromoll@virginia.edu
\end{minipage}
\hfill
\begin{minipage}{2.6in}
\footnotesize {\sc Maria Curie-Sk\l odowska University\\
Department of Mathematics\\
Lublin, Poland\\
E-mail:} lkruk@hektor.umcs.lublin.pl
\end{minipage}
\begin{center}
\begin{minipage}{3in}
\footnotesize {\sc Department of Mathematics \\
California State University, San Marcos \\
San Marcos CA 92096\\
E-mail:} apuha@csusm.edu
\end{minipage}
\end{center}

\end{document}